\documentclass{article}

\usepackage{amssymb,amsmath,mathtools,xcolor,graphicx,xspace,colortbl,ragged2e,rotating} % 
\graphicspath{{../graphics/}{../tcache/}{../gcache/}}
\DeclareGraphicsExtensions{.pdf,.eps,.ps,.png,.jpg,.jpeg}
\newtheorem {theorem}{Theorem}

\newtheorem{definition}[theorem]{Definition}

\newtheorem{lemma}[theorem]{Lemma}

\newenvironment {proof}[1][Proof]{\noindent \textbf {#1.} }{\ \rule {0.5em}{0.5em}}
\begin{document}
\title{Characterizations of pseudo-umbilical submanifolds of locally product Riemannian manifolds}
\author{Ayhan Aksoy}

\maketitle

\begin{abstract}
Compared to totally umbilical submanifolds, studies on pseudo umbilical submanifolds are quite limited. In this paper, pseudo-umbilical submanifolds of locally product Riemannian manifolds are studied.  Necessary and sufficient conditions for a pseudo umbilical submanifold to be invariant, or anti-invariant or semi-invariant submanifold are obtained by using projections arisen  from the product structure of a locally product Riemannian manifold.
\end{abstract}

\newline
\section{Introduction}
\par Submanifolds of almost product or locally product Riemannian manifolds has been introduced by Adati and Bejancu, respectively. Adati studied the invariant, anti-invariant and non-invariant submanifolds of almost product Riemannian manifolds in \cite{adati}. Bejancu defined and studied  the semi-invariant submanifolds of a locally product Riemannian manifold in \cite{bejancu} inspired by his own work about the CR-submanifolds of a Kaehler manifold \cite{bejancu2, bejancu3}. 

\par  After semi-invariant submanifolds of a locally product Riemannian manifold were introduced by Bejancu, there has been several studies on this subject area (for some of these studies, see \cite{matsumoto}, \cite{senlin},  \cite{atceken}, \cite{baymet}, \cite{tastan}). But these studies are quite limited whenever compared to the study of CR-submanifolds of an almost Hermitian (or a Kaehler) manifold.  

\par B. \c{S}ahin consider the problem "Under what conditions a pseudo-umbilical submanifold of a Kaehler manifold is an invariant submanifold, or an anti-invariant submanifold or a CR-submanifold?" in \cite{bay}. In this paper we use a similar approach and give necessary and sufficient conditions for a pseudo-umbilical submanifold of a locally product Riemannian manifold to be an invariant submanifold, or an anti-invariant submanifold or a semi-invariant submanifold in terms of the projections arisen from the product structure.     

\section{Preliminaries}
In this section, we introduce the fundamental notions about the submanifolds of a Riemannian manifold and the notions about the locally product Riemannian manifolds and their submanifolds that will be useful throughout this article. We refer \cite{chen} for the notions about the submanifolds of a Riemannian manifold.      
\par Let $M$ be a submanifold of a Riemannian manifold $(\overline{M}, g)$. The Gauss and Weingarten formulas are given by \begin{equation} \label{eq:1} \overline{\nabla}_{X}Y=\nabla_{X}Y+h(X,Y)\end{equation} and \begin{equation} \label{eq:2} \overline{\nabla}_{X}\xi=-A_{\xi}X+ \nabla_{X}^{\bot}\xi\end{equation} for any $X, Y \in \Gamma \left (TM\right )$ and $\xi \in \Gamma \left (TM^{\bot}\right )$, where $\nabla$ and $\overline{\nabla}$ are Levi-Civita connections of $M$ and $\overline{M}$, respectively, $h$ is the second fundamental form of $M$ and $A_{\xi}$ is the shape operator associated with $\xi$. The second fundamental form and the shape operator are related by \begin{equation} \label{eq:3} g\left( {A}_{\xi}X, Y\right) = g\left( h(X,Y),\xi\right).\end{equation} 
\par Let $M$ be a submanifold of a Riemannian manifold $(\overline{M}, g)$. The mean curvature vector field $H$ of $M$ is defined to be \begin{center}$H={\frac{1}{n}} \underset{i=1}{\overset{n}{\sum}} h(e_{i},e_{i})$\end{center} where $\left\lbrace e_{i}\right\rbrace $ is a local orthonormal basis of $TM$. If $H=0$, then $M$ is called a minimal submanifold. If \begin{equation} \label{eq:4} g\left( h(X,Y), H\right) = g(X,Y) g(H,H) \end{equation} for all $X, Y \in \Gamma \left (TM\right )$, then $M$ is called a pseudo-umbilical submanifold.       
\newline
\par We now recall the notion of locally product Riemannian manifold and some certain type of submanifolds of it.
\par Let $\overline{M}$ be a manifold with a tensor field $F$ of type (1,1) such that \begin{equation} \label{eq:5} F^{2}=I (F \not= \overset{-}{+}I) \end{equation} where $I$ is the identity endomorphism of $\Gamma(T\overline{M})$. In this case the endomorphism $F$ is called an almost product structure on $\overline{M}$. If  $(\overline{M},g)$ is a Riemannian manifold such that for all $X, Y \in \Gamma \left (T \overline{M}\right )$  \begin{equation} \hspace{5 mm} g(X,Y)=g(FX,FY), \end{equation} then $(\overline{M},g)$ is called an almost product Riemannian manifold. Furthermore, if $\overline{\nabla}$ is the Levi-Civita connection associated with the Riemannian metric $g$ and the following identity is satisfied for all  $X, Y \in \Gamma \left (T \overline{M}\right )$ \begin{equation} \label{eq:6} \left( \overline{\nabla}_{X}F\right) Y=0,\end{equation} then $\overline{M}$  is called a locally product Riemannian manifold.
\newline
\par Let $M$ be a Riemannian submanifold isometrically immersed in a locally product Riemannian manifold $(\overline{M},g)$ with the product structure $F$. Then the submanifold $M$ can be characterized according to the action of $F$ on the tangent space $TM$ as follows:
\begin{definition} Let $M$ be a submanifold of a locally product Riemannian manifold $\overline{M}$ with the product structure $F$. If there exist two distributions $D$ and $D^{\bot}$ such that $F(D)=D$,  $F\left(D^{\bot}\right) \subset TM^{\bot}$ and $TM=D \oplus D^{\bot}$, then $M$ is called a semi-invariant submanifold of $\overline{M}$. If $M$ is a semi-invariant submanifold and the dimension of the distribution $D$ (respectively, the dimension of the distribution $D^{\bot})$ is zero, then M becomes an anti-invariant submanifold (respectively, invariant submanifold) of $\overline{M}$. If $M$ is a semi-invariant submanifold such that it is neither an invariant nor an anti-invariant submanifold, then it is called a proper semi-invariant submanifold \cite{baymet}. \end{definition}
\par For all $X \in \Gamma \left (TM\right )$, $FX$ can be decomposed into its tangent part and normal part in the following way: \begin{equation} \label{eq:7} FX=\phi X+\omega X\end{equation} where $\phi X$ and $\omega X$ are the tangent part and normal part of $FX$, respectively. Similarly for each vector field $\xi$ normal to $M$ \begin{equation} \label{eq:8} F\xi=B\xi+C\xi\end{equation} where $B\xi$ and $C\xi$ are the tangent part and normal part of $F\xi$, respectively. 
\newline
\par It follows that $M$ is an invariant submanifold iff $\omega$ is identically zero and $M$ is an anti-invariant submanifold iff $\phi$ is identically zero. The following theorem from (\cite{atceken}, p.28, Theorem 4.1) gives the necessary and sufficient conditions for a submanifold of a locally product Riemannian manifold to be a semi-invariant submanifold. 
\begin{theorem}
	Let $M$ be a submanifold of a locally product Riemannian manifold $\overline{M}$. $M$ is a semi-invariant submanifold iff $\omega \phi=0$ and $rank(\phi)=constant$ \cite{atceken}. 
\end{theorem} 
We also have the following lemmas that will be useful in the following section.
\begin{lemma} Let $M$ be a submanifold of a locally product Riemannian manifold $\overline{M}$ with the product structure $F$. For all $X,Y \in \Gamma(TM)$ the normal part of $\left(\overline{\nabla}_{X} F\right)Y$ is given by
	\begin{equation} \label{eq:10} \left({\nabla}_{X} \omega\right)Y+h(X,\phi Y)=Ch(X,Y)\end{equation}
where $\left({\nabla}_{X} \omega\right)Y=\nabla_{X}^{\bot} \omega Y-\omega \nabla_{X}Y$.
\end{lemma} 
\begin{proof}
For all $X,Y \in \Gamma(TM)$ if we use the equations (\ref{eq:7}), (\ref{eq:1}), (\ref{eq:2}) and (\ref{eq:8}) in the equation (\ref{eq:6}), then we have 
	\begin{Center}
		$\nabla_{X} \phi Y+h(X,\phi Y)-A_{\omega Y}X+\nabla_{X}^{\bot}\omega Y=\phi \nabla_{X}Y+ \omega \nabla_{X}Y+Bh(X,Y)+Ch(X,Y)$.
	\end{Center}
 From the normal part of this equation we have the equation (\ref*{eq:10}).
\end{proof}

\begin{lemma}  Let $M$ be a submanifold of a locally product Riemannian manifold $\overline{M}$ with the product structure $F$. For all $X \in \Gamma(TM)$ and $\xi \in \Gamma(TM^{\bot})$, the normal part of $\left(\overline{\nabla}_{X} F\right)\xi$ is given by 
\begin{equation} \label{eq:12} \left({\nabla}_{X} C\right)\xi=-\omega A_{\xi}X-h(X,B \xi)\end{equation} where $\left ( \nabla _{X}C\right )\xi = \nabla _{X}^{ \bot }C \xi -C \nabla _{X}^{ \bot } \xi$.
\end{lemma}
\begin{proof}
For all $X \in \Gamma(TM)$ and $\xi \in \Gamma(TM^{\bot})$ if we use the equations (\ref{eq:7}), (\ref{eq:1}), (\ref{eq:2}) and (\ref{eq:8}) in the equation (\ref{eq:6}), then we have 
	\begin{Center}
		$\nabla_{X} B \xi+h(X,B \xi)-A_{C \xi}X+\nabla_{X}^{\bot}C \xi=-\phi A_{\xi}X- \omega A_{\xi}X+B\nabla_{X}^{\bot}\xi+C\nabla_{X}^{\bot}\xi$.
	\end{Center}
From the normal part of this equation we have the equation (\ref*{eq:12}).
\end{proof}

\section{Pseudo-umbilical submanifolds of locally product Riemannian manifolds}
In this section we give the main theorems of this article which determine necessary and sufficient conditions for a pseudo-umbilical submanifold of a locally product Riemannian manifold to be an invariant submanifold, an anti-invariant submanifold and a  semi-invariant submanifold, respectively. 
   
\begin{theorem}
Let $M$ be a pseudo-umbilical submanifold of a locally product Riemannian manifold $\overline{M}$. Then for all $X \in \Gamma \left (TM\right )$, the identity $\left ( \nabla _{X}C\right)H=-h(X,BH)$ is satisfied if and only if one of the following is satisfied;
\begin{enumerate}
\item $M$ is a minimal submanifold,
\item $M$ is an invariant submanifold,
\end{enumerate}
where $H$ is the mean curvature vector field and $\left ( \nabla _{X}C\right )H = \nabla _{X}^{ \bot }CH -C \nabla _{X}^{ \bot }H$.
\end{theorem}

\begin{proof}
Since $M$ is a pseudo-umbilical submanifold, then for all $X \in \Gamma(TM)$ we have  ${A}_{H}X={\left \Vert \ H \right \Vert}^2 X$. If we replace $\xi$ with $H$ in the equation (\ref{eq:12}), then we have \begin{Center}  $\left({\nabla}_{X} C\right)H=-{\left \Vert \ H \right \Vert}^2 \omega X-h(X,B H)$. \end{Center}   It follows that $\left ( \nabla _{X}C\right)H=-h(X,BH)$ for all $X \in \Gamma \left (TM\right )$ if and only if either $M$ is a minimal submanifold (i.e. $H=0$) or $M$ is an invariant submanifold (i.e. $\omega=0$)
\end{proof}

Following theorem gives the necessary and sufficient conditions for a submanifold of a locally product Riemannian manifold to be an anti-invariant submanifold.
\newline 
\begin{theorem}
Let $M$ be a pseudo-umbilical submanifold of a locally product Riemannian manifold $(\overline{M},g)$. Then for all $X, Y \in \Gamma \left (TM\right )$, the identity \begin{center} $g\left (\left ( \nabla _{X} \omega \right )Y, H \right)=g \left (Y, A _{CH} X \right )$  \end{center} is satisfied if and only if one of the following is satisfied;
\begin{enumerate}
\item $M$ is a minimal submanifold,
\item $M$ is an anti-invariant submanifold,
\end{enumerate}
where $H$ is the mean curvature vector field and $\left ( \nabla _{X} \omega \right) Y = \nabla _{X}^{ \bot } \omega Y -\omega \nabla _{X} Y$. 
\end{theorem}
\begin{proof}
	If we take the inner product of both sides of the equation (\ref{eq:10}) with $H$, we have \begin{center}
	 $g\left( \left( \nabla_{X} \omega\right) Y,H\right) +g\left( h(X,\phi Y), H\right) =g\left( Ch(X,Y), H\right) $. \end{center} Since $M$ is a pseudo-umbilical submanifold of a product Riemannian manifold, then we have \begin{center} $g\left( \left( \nabla_{X} \omega\right) Y,H\right) + {\left \Vert \ H \right \Vert}^2 g(X, \phi Y)=g\left(Y, A_{CH}X\right)$ \end{center} by the equations (\ref{eq:3}), (\ref{eq:5}) and (\ref{eq:4}). From this last equation it is obvious that $g\left (\left ( \nabla _{X} \omega \right )Y, H \right)=g \left (Y, A _{CH} X \right )$ for all $X,Y \in \Gamma \left (TM\right )$if and only if either $M$ is a minimal submanifold (i.e. $H=0$) or M is an anti-invariant submanifold (i.e. $\phi=0$)
\end{proof}
\newline
\begin{theorem}
	Let $M$ be a pseudo-umbilical submanifold of a locally product Riemannian manifold $\overline{M}$. Then for all $X \in \Gamma \left (TM\right )$, the identity \begin{center} $g\left (\left ( \nabla _{\phi X}C\right )H, CH \right )=-g(h(\phi X, BH), CH) $ \end{center} is satisfied if and only if one of the following is satisfied;
	\begin{enumerate}
		\item $M$ is a minimal submanifold,
		\item $M$ is a semi-invariant submanifold
		\item $\omega \phi X$ is perpendicular to $CH$
	\end{enumerate}
	where $H$ is the mean curvature vector field and $\left ( \nabla _{\phi X}C\right )H= \nabla _{\phi X}^{ \bot } CH - C\left ( \nabla _{\phi X}^{ \bot }H\right )$. 
\end{theorem}
\begin{proof}
	In the equation (\ref{eq:12}) if we replace $X$ with $\phi X$ and replace $\xi$ with $H$, then we have $\left({\nabla}_{\phi X} C\right)H=-\omega A_{H}\phi X-h(\phi X,B H)$ for all $X \in \Gamma(TM)$, where $H$ is the mean curvature vector field. If we take the inner product of both sides of this equation with $CH$ and use the equations (\ref{eq:3}) and (\ref{eq:4}), then we have \begin{Center} $g\left (\left ( \nabla _{\phi X}C\right )H, CH \right )=-{\left \Vert \ H \right \Vert}^2 g\left( \omega \phi X, CH\right) -g(h(\phi X, BH), CH) $. \end{Center} From this last equation it follows that $g\left (\left ( \nabla _{\phi X}C\right )H, CH \right )=-g(h(\phi X, BH), CH) $ for all $X \in \Gamma(TM)$ if and only if one of the following is satisfied:
	\begin{enumerate}
		\item $M$ is a minimal submanifold (i.e. $H=0$),
		\item $M$ is a semi-invariant submanifold (i.e. $\omega \phi=0$),
		\item $\omega \phi X$ is perpendicular to $CH$.
	\end{enumerate}
\end{proof}


\begin{thebibliography}{99}
\bibitem{bejancu2}
A., Bejancu, \emph{CR submanifolds of a Kaehler manifold I}, Proc. Am. Math. Soc. 69(1), 135-142 (1978).
\bibitem{bejancu3}
A., Bejancu, \emph{CR submanifolds of a Kaehler manifold II}, Trans. Am. Math. Soc. 250,333-345 (1979).
\bibitem{bejancu}
A., Bejancu, \emph{Semi-Invariant Submanifolds of locally product Riemannian manifolds}, Ann. Univ. Timisora, s. Math., XXII (1984), 3-11.
\bibitem{chen}
B.Y., Chen, \emph{Differential Geometry of Warped Product Manifolds and Submanifolds}, World Scientific, Hackensack, NJ 2017.
\bibitem{bay}
B., \c{S}ahin, \emph{Characterizations of pseudo-umbilical submanifolds of Kaehler manifolds by means of specific types}, Journal of Geometry, Vol.113 (34), 2022.
\bibitem{baymet}
B., \c{S}ahin and  M., At\c{c}eken, \emph{Semi-invariant submanifolds of Riemannian product manifold}, Balkan Journal of Geometry and Its Applications, Vol.8, No.1, 2003, pp. 91-100.
\bibitem{tastan}
H.M., Tastan and F., Ozdemir, \emph{The geometry of hemi-slant submanifolds of locally product Riemannian manifold}, Turk J. Math. 39 (2015)  
\bibitem{matsumoto}
K., Matsumoto, \emph{On submanifolds of locally product Riemannian manifolds}, TRU Mathematics 18-2, 1982, 145-157.
\bibitem{atceken} 
M., At\c{c}eken, \emph{Geometry of semi-invariant submanifolds of a Riemannian product manifold}, Mathematica Moravica, Vol. 14-1 (2010) 23-34.
\bibitem{adati}
T. Adati, \emph{Submanifolds of an almost product Riemannian manifold}, Kodai Math. J., 4-2 (1981), 327-343.
\bibitem{senlin}
X., Senlin, Yilong, N. \emph{Submanifolds of product Riemannian manifolds}, Acta Mathematica Scientia 2000, 20(B) 213-218.
\end{thebibliography}
\end{document}